\documentclass[12pt]{article}
\newtheorem{lm}{Lemma}[section]

\newtheorem{prop}{Proposition}[section]

\newtheorem{thm}{Theorem}[section]
\newcommand{\qed}{\hspace*{1mm} \rule{2mm}{3mm}\vspace*{0.3cm}}
\newcommand{\be}{\begin{equation}}
\newcommand{\ee}{\end{equation}}
\newcommand{\bea}{\begin{eqnarray}}
\newcommand{\eea}{\end{eqnarray}}
\newcommand{\ba}{\begin{array}}
\newcommand{\ea}{\end{array}}
\newcommand{\beas}{\begin{eqnarray*}}
\newcommand{\eeas}{\end{eqnarray*}}

\newcommand{\ZZ}{{\bf Z}}

\newcommand{\RR}{{\bf R}}
\newcommand{\NN}{{\bf N}}

\newcommand{\DD}{{\cal D}}
\newcommand{\nn}{\nonumber}

\def\etatil{\widetilde{\eta}}

\def\e{\varepsilon}

\def\aihpnl{Ann.\ Inst.\ H.\ Poincar\'e Anal.\ Non Lin\'eaire}

\begin{document}
 \title{An upper bound on the fluctuations of a second class
particle}
 \author{Timo Sepp\"al\"ainen\\
 Department of Mathematics\\
 University of Wisconsin\\
 Madison, WI 53706-1388}

 \thispagestyle{empty}
 \maketitle
 \abstract{This note proves an upper bound for the fluctuations 
of a second-class particle in the totally asymmetric simple
exclusion process. The proof needs a lower tail 
estimate for the last-passage growth
model associated with the exclusion process.
A stronger estimate has been
proved for the corresponding discrete time
 model, but not for the continuous time model
we work with. So we take the needed estimate as a hypothesis.
The process is assumed to be 
 initially   in local 
equilibrium with a slowly varying macroscopic profile. 
The macroscopic initial profile is smooth in a 
neighborhood of the origin where the second-class particle 
starts off, and  the forward characteristic 
from the origin is not a shock.  
 Given these assumptions, the result is that the 
typical fluctuation of the second-class particle is not
of larger order than $n^{2/3}(\log n)^{1/3}$, where $n$ is
the ratio of the macroscopic and microscopic space scales. 
The conjectured correct order should be $n^{2/3}$. 
Landim et al.\ have proved a lower bound of order 
$n^{5/8}$ for  more general
asymmetric exclusion processes in equilibrium. 
Fluctuations in the case of shocks and rarefaction fans
are covered by earlier results of Ferrari--Fontes and 
Ferrari--Kipnis.  
 }

\hbox{} 
 
 \thanks{Research  partially supported by NSF grant DMS-0126775.
 \\[.15cm]
 {\sl Abbreviated title}:  Fluctuation bound  for second class
particle \\[.15cm]
 {\sl AMS (1991) subject classifications}: Primary 60K35; secondary 60F05, 
82C22.}
 \eject

\section{Introduction and result}
We study the motion of a second class particle in a 
 totally asymmetric simple exclusion process on the one-dimensional
integer lattice $\ZZ$. This process describes the evolution of
  indistinguishable particles
that randomly jump to the right on the lattice, one step
at a time. Jumps to already occupied sites are prohibited. 
The state of the process at time $t\ge 0$ 
is the configuration
$\eta(t) =(\eta_i(t):i\in\ZZ)\in\{0,1\}^\ZZ$ 
of occupation numbers, where
$\eta_i(t)=1$ if site $i$ is occupied by a  particle
at time $t$, and $\eta_i(t)=0$ if  site $i$ is vacant 
at time $t$. 
 
 The process is constructed on a probability space
on which are defined the initial particle configuration 
$\eta(0) =(\eta_i(0))_{i\in\ZZ}$, and independently
of $\eta(0)$,  a collection 
 $\{D_i:i\in\ZZ\}$
  of mutually independent rate 1  
Poisson point processes on the time axis
$(0,\infty)$. 
  $D_i$ is the random set of 
 time points (or epochs) when a  jump
from 
site $i$ to site $i+1$ is attempted. Such a jump is
executed at an epoch $t$ of $D_i$ if immediately prior
to time $t$  site $i$ is occupied and site $i+1$ is vacant.
In other words $\eta_i(t-)=1$ and $\eta_{i+1}(t-)=0$, and 
then after the jump   $\eta_i(t)=0$ and $\eta_{i+1}(t)=1$.

The dynamics can be
 represented by the generator $L$  that
acts on bounded
cylinder functions $f$ on the state space $\{0,1\}^\ZZ$:
$$Lf(\eta)=\sum_{i\in\ZZ} \eta_i(1-\eta_{i+1}) 
[f(\eta^{i,i+1})-f(\eta)].
$$
Here $\eta^{i,i+1}$ is the configuration that results from the
 jump of a single particle  from site $i$ to site $i+1$.

The position of a second class particle is defined as follows.
Let $X(0)\in\ZZ$ be a random initial position, and suppose
that initially $\eta_{X(0)}(0)=0$.  
Define a second initial configuration $(\etatil_i(0):i\in\ZZ)$
that differs from $\eta(0)$ only at site $X(0)$: 
$\etatil_{X(0)}(0)=1$ and $\etatil_i(0)=\eta_i(0)$
for $i\ne X(0)$. 
Run the processes  $\eta(t)$ and $\etatil(t)$ 
so that they read the same Poisson jump time processes $\{D_i\}$.
Then there  is always a unique site $X(t)$ at which the two 
processes differ: 
$\etatil_{X(t)}(t)=\eta_{X(t)}(t)+1$ and $\etatil_i(t)=\eta_i(t)$
for $i\ne X(t)$. This defines the position  
$X(t)$  of the second class
particle. 

We refer to the literature for further details of the construction
of these processes. See Chapters III.1--2 in \cite{lig}. 

Now we consider the hydrodynamic limit setting. Assume given
a measurable function $0\le \rho_0(x)\le 1$ on $\RR$. 
Suppose we have
 a sequence $\eta^n=(\eta^n_i)_{i\in\ZZ}$ of random initial configurations
that satisfy, for all finite $a<b$, 
\be
\lim_{n\to\infty}\frac1n \sum_{i=[na]+1}^{[nb]}
\eta^n_i= \int_a^b\rho_0(x)dx
\quad\mbox{in probability.}
\label{hydroass1}
\ee
From this assumption follows a hydrodynamic limit. Let
$\eta^n(t)$ denote the process with initial
configuration $\eta^n$. Then 
for all finite $a<b$ and $0\le t<\infty$, 
\be
\lim_{n\to\infty}\frac1n \sum_{i=[na]+1}^{[nb]}
\eta^n_i(nt)= \int_a^b\rho(x,t)dx
\quad\mbox{in probability.}
\label{hydrolim1}
\ee
The macroscopic profile $\rho(x,t)$ is the unique entropy solution 
of the nonlinear scalar conservation law 
\be
\rho_t+f(\rho)_x=0\,,\quad \rho|_{t=0}=\rho_0, 
\label{conslaw}
\ee
with flux function 
\be
f(\rho)=\rho(1-\rho).
\label{fluxdef}
\ee

To describe the behavior of a second class particle in
this scaling,  
we construct 
the solution $\rho(x,t)$ of (\ref{conslaw}) via 
the Lax-Oleinik formula, and then show how the characteristics of 
(\ref{conslaw}) are defined in this setting.  

Let $g$ be the nonincreasing, nonnegative convex function on $\RR$
defined by 
\be
g(x)=\sup_{0\le\rho\le 1}\{ f(\rho)-x\rho\}
=
\left\{\ba{ll}
-x, &x< -1 \\
(1/4)(1-x)^2, &-1\le x\le 1 \\
0, &x\ge 1.
\ea \right.
\label{fgduality}
\ee

 Define an 
antiderivative $u_0$ of $\rho_0$ by 
\be
u_0(0)=0  \quad\mbox{  and }\quad
u_0(b)-u_0(a)=\int_a^b \rho_0(x)dx\quad\mbox{for all $a<b$.}
\label{defu0}
\ee
For $x\in\RR$ set $u(x,0)=u_0(x)$, and for $t>0$ 
\be
u(x,t)=\sup_{y\in\RR}\left\{u_0(y)-tg\left(\frac{x-y}t\right)\right\}
\label{hopflax}
\ee
The supremum is attained at some $y\in[x-t,x+t]$. The 
function $u$ is uniformly Lipschitz  on $\RR\times[0,\infty)$,
nonincreasing in  $t$ and  nondecreasing in $x$. 
(\ref{hopflax}) is   the {\it Hopf-Lax formula}.
It defines  $u(x,t)$
as the unique viscosity solution of the  Hamilton-Jacobi
equation 
$$u_t+f(u_x)=0\,,\qquad u(x,0)=u_0(x).
 $$

Define the minimal and maximal Hopf-Lax  maximizers in
(\ref{hopflax}) by
$$y^-(x,t)=\inf\left\{ y\ge x-t: u(x,t)=u_0(y)-t
g\left(\frac{x-y}{t}\right) \right\}
$$
and 
$$y^+(x,t)=\sup\left\{ y\le x+t: u(x,t)=u_0(y)-t
g\left(\frac{x-y}{t}\right)  \right\}. 
$$
The entropy solution
of (\ref{conslaw}) 
 is defined by
the {\it Lax-Oleinik formula:}
\be
\rho(x,t) =
-g'\left(\frac{x-y^\pm(x,t)}{t}\right).  
\label{defrho}
\ee
This definition makes sense a.e.\ because for a fixed $t$, 
$y^-(x,t)=y^+(x,t)$ for all but countably many $x$. 
  The derivative 
$u_x(x,t)$ exists and equals $\rho(x,t)$ for all $(x,t)$ such 
that $y^-(x,t)=y^+(x,t)$. 
A point $(x,t)$ for $t>0$ is a {\it shock} if $y^-(x,t)\ne y^+(x,t)$.
Equivalently, $\rho$ is not continuous at $(x,t)$. 

The minimal and maximal forward characteristics are
defined for $b\in\RR$, $t>0$, as
\be
w^-(b,t)=\inf\{ x: y^+(x,t)\ge b\}
\ \mbox{ and }\ 
w^+(b,t)=\sup\{  x: y^-(x,t)\le b\}. 
\label{defw}
\ee
 The forward characteristics $w^\pm(b,t)$ are Filippov solutions 
of the ordinary differential equation 
$dx/dt=f'(\rho(x,t))$, $x(0)=b$. If $\rho_0$ is continuous
at $b$, the forward characteristic is unique, in other words
$w(b,t)=w^\pm(b,t)$. 

Now return to the hydrodynamic limit setting where a sequence
of processes $\eta^n(t)$ is assumed to satisfy  (\ref{hydroass1}). 
Let $X_n(t)$ be the position of a second class particle in
 process  $\eta^n(t)$. Assume the initial location is always the
origin: $X_n(0)=0$ a.s. 
Assume also that the  forward characteristic from
$0$  is unique, so $w(0,t)=w^\pm(0,t)$. 
Then  we have the law
of large numbers
\be
\lim_{n\to\infty}\frac1n X_n(nt)= w(0,t)
\qquad\mbox{in probability}. 
\label{Xntlimit}
\ee
Original references for the hydrodynamic limits above are
 \cite{rez1, sepK}, and for the second class particle
limit  \cite{fer,  rez2,  septra}.
See 
  \cite{ kiplan, spo} for  general treatments of the 
macroscopic behavior of the exclusion process, and 
\cite{eva} for the basic p.d.e.\ theory used above. 

The result of our note is on the fluctuations from the 
limit (\ref{Xntlimit}). First a brief mention of known 
results. Let $\nu_{\lambda,\rho}$ denote the product measure
on the state space $\{0,1\}^\ZZ$ under which particles are
present with density $\lambda$ ($\rho$) to the left (right)
of the 
origin. Put a second class particle initially 
at the origin. Now there is only a single process,
not a sequence.  Consider two cases. 

(i) 
 If $\lambda<\rho$ then the characteristic
from the origin, $w(0,t)=t(1-\lambda-\rho)$, is a shock. 
In this case $X(t)$ has diffusive fluctuations. 
It satisfies
a central limit theorem with explicitly known variance
 in the scale $\sqrt{t}$. This result is from
Ferrari and Fontes \cite{ferfon}, and is covered in
Liggett's monograph \cite{lig}. 

(ii)  If $\lambda>\rho$ then the characteristic
from the origin is not unique, and we have 
 $w^-(0,t)=t(1-2\lambda)$,  $w^+(0,t)=t(1-2\rho)$.
$X(t)$ has macroscopic fluctuations: $t^{-1}X(t)$ converges
weakly to a uniform distribution on $[w^-(0,t), w^+(0,t)]$. 
This result is due to Ferrari and Kipnis \cite{ferkip}. 

In particular, there is currently no definitive result in the 
equilibrium case $\lambda=\rho$. In this case it is 
conjectured that the typical fluctuations are of order
$t^{2/3}$. Recent results on the $t^{1/3}$ fluctuations for
growth models provide corroboration for this conjecture, 
see \cite{praspo}. 
Landim et al. \cite{lqsy} proved a $t^{5/4}$  lower bound on the variance
in the following weak sense: there exists a constant $C>0$ 
such that for small enough $\lambda>0$, 
$$
\int_0^\infty e^{-\lambda t}\mbox{ Var}(X(t))dt \ge C\lambda^{-9/4}.
$$ 
A large deviation bound 
for the second class particle 
in the equilibrium situation was proved in \cite{sepset},
and used to prove central limit theorems for additive functionals
of the exclusion process. 

We prove an upper bound on the size of the typical
 fluctuation in the limit (\ref{Xntlimit}). Make the 
following assumptions on the macroscopic profile around
the origin: 
\be
\ba{ll}
&\mbox{$\rho_0$ is continuously differentiable
in a neighborhood of the origin, and   }\\
&\mbox{$0<\rho_0(0)<1$. Suppose
 $t>0$ is such that $\rho_0'(0)<1/(2t)$.}
\ea
\label{assrho0}
\ee
Smoothness around the origin forces the forward characteristic
to be unique, so we can denote this value by $w(0,t)=w^\pm(0,t)$.
The assumptions $0<\rho_0(0)<1$ and $\rho_0'(0)<1/(2t)$
are for technical convenience: the first implies that
$w(0,t)$ lies in the interior of $(-t,t)$, and the second
that the mapping from  $x\to y^\pm(x,t)$ is Lipschitz both ways
for $x$ close enough to $w(0,t)$. 

Secondly, we need to assume that $w(0,t)$ is not a shock:
\be
y^-(w(0,t),t)=y^+(w(0,t),t).
\label{notshock}
\ee

The limit (\ref{Xntlimit}) required only assumption
(\ref{hydroass1}). We need to 
strengthen this in order to have sharper control over initial
fluctuations. We assume that the initial particle configurations are
in local equilibrium with macroscopic profile $\rho_0$: 
\be
\ba{ll}
&\mbox{For each fixed $n$, the random 
variables $(\eta^n_i)_{i\ne 0}$ are mutually  }\\
&\mbox{independent, and for each $i$, 
$\displaystyle  E[\eta^n_i]=n\int_{(i-1)/n}^{i/n}\rho_0(x)dx$.} 
\ea
\label{loceqass}
\ee

The second class particle starts at the origin, so the origin is
left empty in the initial exclusion configurations [$\eta^n_0(0)=0$].

The proof of our fluctuation bound relies on three 
ingredients. (i) A bound on the typical lower tail
fluctuation of the last-passage growth model associated
with the asymmetric exclusion process. (ii) A sharper 
bound on the upper tail deviations of the growth model.
(iii) A variational representation for the location of the 
second class particle. 

We face a dilemma: ingredients (ii) and (iii),
 the upper tail bound and the variational
representation, have been proved for the totally asymmetric 
exclusion process in continuous time 
considered here \cite{sepmprf, septra}.
But point (i),  the lower tail bound, has presently been proved
only for discrete time exclusion with geometric 
waiting times.   The 
estimate we need  is  just a little more than what follows
from 
Johansson's distributional limit \cite[eqn.\ (1.22)]{joh} for 
continuous time exclusion.  
 Baik et al.\ \cite{bdm} have proved
a much stronger bound for discrete time exclusion,
and there should be no doubt that an analogous estimate 
is true for continuous time too. 
For this reason we feel comfortable in taking the 
 lower tail estimate we need as an extra 
hypothesis.

Consider the following last-passage growth model. 
 Let $\NN=\{1,2,3,\ldots\}$ be the set of
natural numbers, and  
  let $\{u_{i,j}:(i,j)\in\NN^2\}$ be i.i.d.\ exponential
mean 1 random variables. Set 
\be
H(M,N)=\max_\sigma\sum_{(i,j)\in\sigma}u_{i,j},
\label{tHmax}
\ee
where the maximum is over lattice paths 
$$\sigma
=\{ (1,1)= (i_1,j_1), (i_2,j_2),
\ldots,(i_{M+N-1},j_{M+N-1})=(M,N)\}$$
 in $\NN^2$
 that take  
only up-right steps:
$$
\mbox{$(i_{m+1},j_{m+1})-(i_m,j_m)= (0,1)$ or $(1,0)$
for each $m$.}
$$
  Johansson \cite[Theorem 1.6]{joh}
proved a  distributional
limit   for $H([\alpha n],[\beta n])$ as $n\to\infty$,
$(\alpha,\beta)\in\RR^2_+$. What we need for our proof is the following
estimate. 

\medskip

{\bf Hypothesis H.} Suppose $\alpha_n\to\alpha>0$
and $\beta_n\to\beta>0$ are convergent sequences
of positive numbers  such that, for a constant $B<\infty$,
\be
|\alpha_n-\alpha| + |\beta_n-\beta| \le Bn^{-1/3}(\log n)^{1/3}
\qquad\mbox{for all $n$.}
\label{hypoHab}
\ee
Let $\e>0$. Then, if $C$ is fixed large enough, 
\be
P\left\{ H([n\alpha_n],[n\beta_n]) <n(\sqrt{\alpha_n}+\sqrt{\beta_n}\,)^2
-Cn^{1/3}\,\right\} \le \e
\qquad\mbox{for all $n$.}
\label{hypoH}
\ee

\medskip

The connection of $H(M,N)$
 to exclusion is this: start the exclusion process
with all sites to the left of site 1 occupied, and all sites
to the right of site 0 empty. Then the first time when 
there are $j$ particles to the right of site $i$ is
distributed as $H(i+j,j)$ for $j>(-i)\vee 0$. 
  Rost  \cite{ros} treated
the hydrodynamic limit of this particular 
exclusion  process without
  the last-passage formulation. 
There seems to be no first paper to cite
as the original source of the  last-passage connection. 
This author began using the last-passage representation
  in 1995, by which time it was certainly known
by a number of people. In one  conversation  A. Gandolfi
was credited with making this  observation. 

An alternative course 
for this paper would
be to  reprove ingredients (ii) and (iii) for discrete
time exclusion, following 
 \cite{sepmprf, septra}. Then in place of ({\ref{hypoH})
  we could   use the  Baik et al.\ \cite{bdm} 
estimate for the discrete time growth model.
Our theorem would then be for  
 discrete time exclusion. We chose the present course 
since continuous time exclusion is the process that most 
people prefer to work with.  And also to avoid  extra work, 
since in any case our approach does not quite get the optimal order.

We now state the theorem.

\begin{thm} Assume {\rm (\ref{assrho0})}, {\rm (\ref{notshock})},  and 
 {\rm (\ref{loceqass})}, and assume 
Hypothesis H. Let  $\eta^n(t)$ be the totally asymmetric
simple exclusion process with initial configuration $\eta^n$. 
 Let $X_n(t)$ be the location of a second class 
particle in the process $\eta_n(t)$, started at the origin
$X_n(0)=0$. 
 Let $\e>0$. Then there exists a constant
$b<\infty$ such that, for all $n$, 
$$
P\left\{ \,\left| X_n(nt)-nw(0,t)\right|\ge bn^{2/3}(\log n)^{1/3}\right\}
\le \e.
$$
\label{ubound}
\end{thm}

 The rest of the note contains the proof. We start with
properties of the Hopf-Lax formula. Then we explain the 
variational coupling representation of the  totally asymmetric
exclusion process. 
This representation is used in 
conjunction with two probability estimates,   one
  from   \cite{sepmprf},
the other 
from hypothesis (\ref{hypoH}). 

\section{Properties of the macroscopic profile}

Recall the connections (\ref{defu0}) and (\ref{defrho}).
Abbreviate $r=w(0,t)$ throughout the proof. Let 
$$
I(x,t)=\{ y\in\RR: u(x,t)=u_0(y)-t
g((x-y)/{t})\}$$
 denote the set of maximizers in the Hopf-Lax 
formula (\ref{hopflax}). 

\begin{lm} Suppose 
 $u_0$ is twice continuously differentiable
in a neighborhood of $0$. Assume that  $y^\pm(r,t)=0$
   and 
$u_0''(0)<1/(2t)$. Let $a_0>0$.  Then there exist
$\delta_1>0$ and  $0<c_0<\infty$  such that the following
holds: if
$x\in [r-\delta_1, r+\delta_1]$,  
 $y\in I(x,t)$, and  $ \eta\in[-a_0, a_0]$,  then 
\be
u_0(y)-tg\left(\frac{x-y}t\right) \ge u_0(\eta)-
tg\left(\frac{x-\eta}t\right) + c_0(\eta-y)^2.
\label{quadratic1}
\ee
Furthermore, if $x_0\in [r-\delta_1, x)$ is distinct from
$x\in [r-\delta_1, r+\delta_1]$, and  $y_0\in I(x_0,t)$, then 
 \be  
 c_0 (x-x_0)\le y-y_0 \le c_0^{-1} (x-x_0).  
\label{yplusass}
\ee
\end{lm}

{\it Proof.} We first prove 
(\ref{yplusass}). 
By the assumption
$y^\pm(r,t)=0$, $I(x,t)\to \{0\}$ as $x\to r$. 
So for $x$ close enough to $r$ and $y\in I(x,t)$,  
\be
 u'_0(y)=-g'((x-y)/t)=(t-x+y)/(2t),
\label{derivatives1}
\ee
and similarly for $x_0$ and $y_0$. This shows $y\ne y_0$ if $x\ne x_0$. 
Hence 
$$
\frac{x-x_0}{y-y_0}= 1-2t\cdot\frac {u'_0(y)-u'_0(y_0)}{y-y_0}
=1-2t u_0''(\theta)
$$
by the mean value theorem, where $\theta\in(y_0,y)$. 
If $y,y_0$ are close enough to $0$, $1-2t u_0''(\theta)$ is bounded
and bounded away from 0. This can be achieved by taking 
$x,x_0 \in [r-\delta_1, r+\delta_1]$ for small enough $\delta_1$. 
This proves 
(\ref{yplusass}).

By the assumptions, we can choose $\alpha, \delta_2>0$ so that 
$u_0$ is $C^2$ and $u''_0<1/(2t)-\alpha$
 on $(-3\delta_2,3\delta_2)$. Shrink the $\delta_1>0$ chosen in the 
previous paragraph so that 
$I(x,t)\subseteq(-\delta_2,\delta_2)$ for all $x\in[r-\delta_1,r+\delta_1]$. 
By continuity and compactness, there exists $\delta_3>0$ such that
$$
u_0(y)-tg\left(\frac{x-y}t\right) \ge u_0(\eta)-
tg\left(\frac{x-\eta}t\right) + \delta_3
$$
for all $x\in[r-\delta_1,r+\delta_1]$, $y\in I(x,t)$ and
 $2\delta_2\le |\eta|\le  a_0$. Thus for 
these $\eta$, (\ref{quadratic1}) holds with 
$c_0=\delta_3/(a_0+\delta_2)^2$. 

Keeping still  $x\in[r-\delta_1,r+\delta_1]$ and  $y\in I(x,t)$
so that $y\in (-\delta_2,\delta_2)$, 
for $\eta\in(-2\delta_2,2\delta_2)$ use Taylor's theorem with the 
Lagrange form of the remainder term \cite[p.\ 195]{str}:
\beas
u_0(\eta)-
tg\left(\frac{x-\eta}t\right) &=&  u_0(y)-tg\left(\frac{x-y}t\right)
+\frac12 (\eta-y)^2 \left\{u_0''(y+ \theta(\eta-y))-\frac1{2t}\right\}\\
&\le&
 u_0(y)-tg\left(\frac{x-y}t\right) -\frac{\alpha}2 (\eta-y)^2.
 \eeas
Here $\theta\in(0,1)$ depends on $\eta$ and $y$, but the 
upper bound on $u_0''$ works in all cases because 
$y+ \theta(\eta-y)\in (-2\delta_2,2\delta_2)$. 
\qed

\section{The variational coupling}

We summarize briefly the variational coupling representation of 
the process and the second class particle \cite{sepmprf, sepK, septra}.
 Now the exclusion process is constructed in terms of a process
  $z(t)= (z_i(t))_{i\in{\ZZ}}$
  of labeled particles that move on  ${\ZZ}$
subject to the constraint
\be
0\le z_{i+1}(t)-z_i(t)\le 1\quad\mbox{for all $i\in\ZZ$ and $t\ge 0$.}
\label{tz1}
\ee
   In the graphical construction,
  $z_i$ attempts to jump one step to the {\it left} at
epochs of $D_i$. A jump is suppressed if it leads to a
violation of (\ref{tz1}).  We start the process so that
$z_0(0)=0$. The connection between the exclusion $\eta(t)$ and 
the process $z(t)$ is that 
\be
\eta_i(t)=z_i(t)-z_{i-1}(t)\,.
\label{tetaz}
\ee
 This equation is used both ways: given the initial 
configuration $\eta(0)$, define initial $z(0)$ by $z_0(0)=0$ 
and by (\ref{tetaz}).
Construct the process $z(t)$ by the graphical representation. 
And then  define $\eta(t)$ by  (\ref{tetaz}). 

Construct
  a family $\{w^k(t):k\in{\ZZ}\}$ of auxiliary processes. 
 The initial configuration $w^k(0)$ depends on the initial
position $z_k(0)$ through a global shift:
\bea
w^k_i(0)=\left\{ \begin{array}{ll}
z_k(0)\,, &i\ge 0\\
z_k(0)+i\,, &i<0.
\end{array}
\right.
\label{twinit}
\eea
The processes $\{w^k(t)\}$ are coupled to each other
and to $z(t)$ through the
Poisson processes $\{D_i\}$, so that particle $w^k_i$  attempts jumps
to the left at epochs of   $D_{i+k}$. 
 The key  variational coupling property says that 
for all $i\in{\ZZ}$ and $t\ge 0$,
$$
z_i(t)=\sup_{k\in{\ZZ}} w^k_{i-k}(t) \quad\mbox{ a.s.}
$$

 It is convenient to
 decompose $w^k(t)$ into a sum of the initial position
defined by (\ref{twinit}) and the increment determined
by  the Poisson  processes.
Define a family of processes $\{\xi^k(t)\}$ by
$$\mbox{$\xi^k_i(t)=z_k(0)-w^k_i(t)$ for $i\in\ZZ$, $t\ge 0$.}$$
The process
$\xi^k(t)$ does not depend on $z_k(0)$, and depends on
the superscript $k$ only through a translation of the
$i$-index of the Poisson processes
$\{D_i\}$.
  Initially
\bea
\xi^k_i(0)=\left\{ \begin{array}{ll}
0\,, &i\ge 0\\
-i\,, &i<0.
\end{array}
\right.
\label{txiinit}
\eea
 We can think of $\xi^k$ as a growth model on the upper half plane,
so that $\xi^k_i$ gives the height of the interface above
site $i$.
Its dynamics are specified by saying 
 that  $\xi^k_i$
advances one step up at each epoch of $D_{k+i}$, 
provided   these inequalities
are preserved:
\be
\xi^k_i(t)\le \xi^k_{i-1}(t)\qquad\mbox{and}\qquad
\xi^k_i(t)\le \xi^k_{i+1}(t)+1.
\label{txiineq}
\ee
The connection of $\xi^k$ with the exclusion process is this: start the 
exclusion process so that all sites from $k$ to the left are occupied,
and all sites from $k+1$ to the right are vacant. Then 
$\xi^k_i(t)$ is the number of particles to the right of site
$k+i$ at time $t$. 

In terms of $\xi$, the variational coupling can be expressed as
 \be
z_i(t)= \sup_{k\in\ZZ}\{ z_k(0)-\xi^k_{i-k}(t) \}
\label{zvar1}
\ee

Now suppose $\eta_{X(0)}(0)=0$ and  put a second class
particle initially at the (possibly random) location $X(0)$. 
 Then the later location of the  second class particle 
satisfies 
\be
X(t)=\inf\{ i\in{\ZZ}: \mbox{$z_i(t)=z_k(0)-\xi^k_{i-k}(t) $ for
some $k\ge X(0)$}\}.
\label{Xtvar}
\ee
 
We conclude this overview with  two 
 monotonicity properties. 
 First, the coupling of the $\xi^k$ processes through
common Poisson clocks gives us this inequality:
\be
\xi^k_{i-k}(t)\le \xi^l_{i-l}(t)\ \ \mbox{for all $i$, if $k\le l$.}
\label{xikl}
\ee
Second, the
variational coupling has this property.

\begin{lm} The following statement holds almost surely.
 Suppose $z_i(t)=w^k_{i-k}(t)$, and $j>i$. Then 
there exists an index $m$ such that $m\ge k$ and 
$z_j(t)=w^m_{j-m}(t)$. 
\label{ikjmlm}
\end{lm}

{\it Proof.} First check that the statement is true at $t=0$,
by the definition of the family $\{w^k\}$ and the restriction
$0\le z_{i+1}-z_i\le 1$. 

To prove it up to time $t_0$, fix
indices $i_0<<0<<i_1$ so that the Poisson jump time
processes $\DD_{l}$ are empty up to time $t_0$ for 
$l\in\{i_0, i_0+1, i_1, i_1+1\}$. Since the positions
$z_l$ and $w^n_{l-n}$ for $l\in\{i_0, i_0+1, i_1, i_1+1\}$ and 
$n\in\ZZ$ do not change up to time $t_0$, the statement
of the lemma holds for $(i,j)=(i_0, i_0+1)$ and $( i_1, i_1+1)$
and $0\le t\le t_0$. Now prove it simultaneously for
pairs $(i,j)$, $i_0<i<j\le i_1$, by induction on the 
finitely many jump times in $\cup_{i_0<i<i_1}\DD_i\cap(0,t_0]$.

Since the indices $i_0<<0<<i_1$ can be chosen arbitrarily far 
away from the origin, this way the lemma is proved for all  $i<j$
   up to time $t_0$. 
\qed

\section{Fluctuation bounds for $\xi$}

 Here is the lower tail bound for  $\xi$ that we need. 
Recall the last-passage model $H(M,N)$
 discussed in the introduction.  Its precise connection with
$\xi$ is 
\be
P\{\xi_i(t)< j\}=  P\{ H(i+j,j)> t\}.
\label{xiH}
\ee
 
\begin{prop} Let $t>0$ and $\e>0$. Then 
there is a finite constant $C>0$ such that, for all 
$x\in[-t+\e, t-\e]$, 
 all small enough $h>0$ and all $n$, 
$$
P\left\{ \xi_{[nx]}(nt)
\le  ntg\left(\frac{x}{t}\right)-2nh\right\} 
\le \exp\left\{-n\left(\frac{4\sqrt{2}}{3}\cdot \frac{\sqrt{t-x\,}}{t+x}
\cdot h^{3/2} - Ch^{5/2}\right)\right\}.
$$
\label{ldboundprop}
\end{prop}

{\it Proof.}
 Accoding to Theorems 4.2 and 
4.4 in \cite{sepmprf}, 
\be
\sup_n \frac1n\log P\{ H([nr], [nw])>nt\} 
=\lim_{n\to\infty}  \frac1n\log P\{ H([nr], [nw])>nt\} 
=-\Psi_{w,t}(r)
\label{ldboundH}
\ee
where the rate function is defined by 
\bea
\Psi_{w,t}(r) &=&\sqrt{ (t-r-w)^2-4rw\;}\nn\\
&&\qquad-2r\cosh^{-1} \left( \frac{t+r-w}{2\sqrt{tr\,}}\right)
-2w\cosh^{-1} \left( \frac{t+w-r}{2\sqrt{tw\,}}\right)
\label{Psidef}
\eea
for values $r,w,t\ge 0$ that satisfy $\sqrt w + \sqrt r\le \sqrt t$. 
Note that the limiting value is 
\be
(\sqrt w + \sqrt r  \,)^2=\lim_{n\to\infty}  \frac1n H([nw], [nr]). 
\label{Hlim}
\ee

In case the reader wishes to compare with
the source  \cite{sepmprf}, note
that the  last-passage model is indexed slightly differently
there. 
  Our variable $w$ appears 
as  a negative variable $x$ in  \cite{sepmprf}. 
 Estimate (\ref{ldboundH}) was in principle  
 covered later by Johansson's \cite{joh} bigger results, 
but the rate function is not explicitly available in 
 \cite{joh}.

 Expanding $\Psi_{w,t}(\cdot)$ close to the   value
$u=(\sqrt t - \sqrt w  \,)^2$ gives 
\be
\Psi_{w,t}(u-h)=\frac43\cdot \frac{w^{1/4}}{t^{1/4}u^{1/2}}\cdot h^{3/2}
+O(h^{5/2})
\label{Psiexp1}
\ee
where the $O$-term is uniform over $w$ and $u$  bounded 
away from $0$, and $h>0$ is small enough. 
Finally for large  enough $n$
so that $nh>1$, 
\beas 
&&P\left\{ \xi_{[nx]}(nt)
\le ntg\left({x}/{t}\right)-2nh\right\}\\
&\le& P\left\{ \xi_{[nx]}(nt)
< [ntg\left({x}/{t}\right)-nh]\right\}\\
&\le& P\left\{ H([nx+ntg\left({x}/{t}\right)-nh],
[ntg\left({x}/{t}\right)]) >nt\right\}\\
&\le& \exp\left\{ -n\Psi_{w,t}(u-h)\right\},
\eeas
where  $w=tg(x/t)=(t-x)^2/(4t)$ and $u=(\sqrt t - \sqrt w  \,)^2
=x+tg(x/t)$. Bounding
$w$ and $u$ away from $0$ is the same as bounding $x\in(-t,t)$ 
away from both $t$ and $-t$. 
\qed

\section{Proof of Theorem \ref{ubound}}

Set $u_n=n^{2/3}(\log n)^{1/3}$, and still $r=w(0,t)$. 
We prove: given $\e>0$, there exists a constant 
$b<\infty$ such that 
\be
P\{ X_n(nt)\le nr -bu_n\}\le \e
\label{upperest1}
\ee
for all $n$. The corresponding bound on the right,
$$
P\{ X_n(nt)\ge nr +bu_n\}\le \e
$$
follows by an argument so similar it is 
not worth repeating.

Let $\e_0=\e/4$.  The assumption $0<u_0'(0)=\rho_0(0)<1$ entails
that 
\be
-t+3\alpha_1\le r\le t-3\alpha_1
\label{ralpha}
\ee
 for some $\alpha_1>0$ [by (\ref{derivatives1}) applied to 
$(x,y)=(r,0)$].

\begin{lm}
For sufficiently
large  $n$, 
 \be
P\left\{ \mbox{ $z^n_{[nr-bu_n]}(nt)= z^n_k-\xi^k_{[nr-bu_n]-k}(nt)$ for some 
 $k\ge \alpha_1n$}  \right\}\le \e_0. 
\label{avar0}
\ee
\end{lm}

{\it Proof.} By Lemma \ref{ikjmlm}, if 
 $z^n_{[nr-bu_n]}(nt)= z^n_k-\xi^k_{[nr-bu_n]-k}(nt)$ for some 
 $k\ge \alpha_1n$, then also
 $z^n_{[nr]}(nt)= z_m-\xi^m_{[nr]-m}(nt)$ for some 
 $m\ge \alpha_1n$. By passing to the hydrodynamic
 limit and by compactness,
this would imply the existence of $\eta\ge \alpha_1$ 
such that 
$u(r,t)=u_0(\eta)-tg((r-\eta)/t)$. See 
\cite[Prop.\ 5.1]{septra} for more details on this type of an
argument. But now we have a contradiction, because by
assumption $0$ is the unique maximizer for $u(r,t)$ in the 
Hopf-Lax formula (\ref{hopflax}). 
\qed

 So now,  except on an event of probability  less than $\e_0$, 
the inequality 
\be
X_n(nt)\le nr -bu_n
\label{goalineq}
\ee
 implies that
 \be
z^n_{[nr-bu_n]}(nt)= z^n_k-\xi^k_{[nr-bu_n]-k}(nt) 
\label{avar1}
\ee
 holds 
for some $0\le k\le \alpha_1n$, if $n$ is large enough. 
The follows from combining variational formula (\ref{Xtvar})
with Lemmas \ref{ikjmlm} and \ref{avar0}. 

Let $y_n\in I(r-bu_nn^{-1}, t)$ be a Hopf-Lax maximizer for
the point
$(r-bu_nn^{-1}, t)$.
  Note that $y_n<0$, and $y_n\nearrow 0$ as 
$n\nearrow\infty$. 
Together with (\ref{avar1}),  the variational equation (\ref{zvar1})
implies  that 
 \be
 z_{[ny_n]}^n-\xi^{[ny_n]}_{[nr-bu_n]-[ny_n]}(nt) \le 
z^n_k-\xi^k_{[nr-bu_n]-k}(nt)
\label{avar2}
\ee
for some $0\le k\le \alpha_1n$. 
For each $n$, set $k_j=[jn^{1/3}]$ for $0\le j\le \alpha_1n^{2/3}$.  
For $k_j\le k\le k_{j+1}$,  we have $z^n_k\le z^n_{k_{j+1}}$,
and   (\ref{xikl}) gives 
$$
\xi^{k_j}_{[nr-bu_n]-k_j}(nt)
\le \xi^k_{[nr-bu_n]-k}(nt).  
$$
From (\ref{avar2}) we then get an intermediate conclusion that,
on the event (\ref{goalineq}), 
$$
\xi^{k_j}_{[nr-bu_n]-k_j}(nt)
\le
\xi^{[ny_n]}_{[nr-bu_n]-[ny_n]}(nt) +  z_{k_{j+1}}- z_{[ny_n]}
$$
for some $0\le j\le \alpha_1n^{2/3}$, except for an event of probability less
than $\e_0$. 

Now we estimate the random variables on the right-hand side 
of the last inequality. 
By (\ref{ralpha}) and (\ref{yplusass}), for large  $n$ the point 
$r-bu_nn^{-1}-y_n$ lies in $[-t+\alpha_1, t-\alpha_1]$. Having this
point in the interior of $(-t,t)$ permits us to apply
the hypothesis  made in the Introduction. Note that 
$bu_nn^{-1}$ is of order $O(n^{-1/3}(\log n)^{1/3})$,
and so is $y_n$ by (\ref{yplusass}). Thus  (\ref{hypoHab})
is satisfied. Then 
by  (\ref{hypoH}) and (\ref{xiH}),  
we may fix a constant $a_1<\infty$
such that, for all $n$, 
$$
P\left\{
\xi^{[ny_n]}_{[nr-bu_n]-[ny_n]}(nt) \le 
ntg(t^{-1}(r-bu_nn^{-1}-y_n))+a_1n^{1/3}
\right\}
\ge 1-\e_0.
$$
 For the initial state, if 
we increase $a_1$ sufficiently, then 
\beas
&&P\left\{ \,
 |z^n_{k_j}-z^n_{[ny_n]}-nu_0(k_j/n)+nu_0(y_n)|\le a_1(k_j-ny_n)^{1/2}
(\log n)^{1/2}\right. \\
&& \left. \qquad\qquad\qquad \qquad\qquad\qquad\qquad\qquad\qquad
\qquad\qquad  {}^{} \mbox{ for all $0\le j\le \alpha_1n^{2/3}$}
\right\}  \ge 1-\e_0
\eeas
for all $n$. This is proved by an 
exponential Chebychev's inequality, with separate arguments 
for the two sides.

Now outside an event of probability  less than $3\e_0$, 
  it  follows from (\ref{goalineq}) that 
\be
\ba {rl}
\xi^{k_j}_{[nr-bu_n]-k_j}(nt)
\le& ntg(t^{-1}(r-bu_nn^{-1}-y_n)) 
+n(u_0(k_{j+1}/n)-u_0(y_n))\\
& \qquad  +a_1n^{1/3}  + a_1(k_{j+1}-ny_n)^{1/2}
(\log n)^{1/2} 
\ea
\label{avar3}
\ee
for some $0\le j< \alpha_1n^{2/3}$. 

 Set $x_j=r-bu_nn^{-1}-k_{j}n^{-1}$.
 Use (\ref{quadratic1}) to write  
\beas
 tg\left(\frac{r}{t}-\frac{bu_n}{nt}-\frac{y_n}{t}\right)
+u_0(k_{j+1}/n)-u_0(y_n)&\le& tg\left(
\frac{r}{t}-\frac{bu_n}{nt}-\frac{k_{j+1}}{nt}\right)
  -
a_2\left(\frac{k_{j+1}}{n}-y_n\right)^2\\
&=& tg(x_{j+1}/t) -
a_2\left(\frac{k_{j+1}}{n}-y_n\right)^2,
\eeas
for a constant $a_2$. 
 Then (\ref{avar3}) gives
 \be
\ba {rl}
\xi^{k_j}_{[nx_j]}(nt)
\le& ntg(x_j/t) +a_1n^{1/3} 
+ a_1(k_{j+1}-ny_n)^{1/2}(\log n)^{1/2}  \\
& \quad -a_2n(k_{j+1}/n-y_n)^2 + 
ntg(x_{j+1}/t)-ntg(x_j/t). 
 \ea
\label{avar4}
\ee
 By the Lipschitz property of $g$, 
$$
ntg(x_{j+1}/t)-ntg(x_j/t)\le a_3n^{1/3},
$$
for a constant $a_3<\infty$. 
Use  (\ref{yplusass}) in the form
$$y_n<-\delta_1 bu_nn^{-1}=-\delta_1 bn^{-1/3}(\log n)^{1/3}.$$ 
Substitute in $k_j=[jn^{1/3}]$ and 
 $u_n=n^{2/3}(\log n)^{1/3}$, and  use the inequalities
$(p+q)^{1/2}\le p^{1/2}+q^{1/2}$ and
$(p+q)^{2}\ge p^{2}+q^{2}$ valid for $p,q\ge 0$.  
 After these steps,  the right-hand side of (\ref{avar4}) is bounded above by
\beas
&&ntg(x_j/t) +(a_1+a_3)n^{1/3} 
+ a_1(j+1)^{1/2} n^{1/6}(\log n)^{1/2}   +a_1n^{1/2}|y_n|^{1/2}(\log n)^{1/2} 
 \\
&& \qquad  -a_2(j+1)^2n^{-1/3} -a_2ny_n^2  \\
&=&ntg(x_j/t)  +\left\{ (a_1+a_3)n^{1/3} - 
 (a_2b^2 \delta_1^2/4)   n^{1/3}(\log n)^{2/3} \right\} \\
&& \quad +\left\{  a_1(j+1)^{1/2} n^{1/6}(\log n)^{1/2} 
-   a_2(j+1)^2n^{-1/3} - (a_2b^2 \delta_1^2/4)   n^{1/3}(\log n)^{2/3}
 \right\} \\
&& \quad +\left\{ a_1n^{1/2}|y_n|^{1/2}(\log n)^{1/2} 
- (a_2/4)ny_n^2  \right\}\\
&& \quad -  (a_2b^2 \delta_1^2/4)   n^{1/3}(\log n)^{2/3} 
\eeas
By fixing $b$ large enough relative to  the 
other  
constants $\delta_1$ and $a_1$--$a_3$, the expressions in braces can be 
made nonpositive for all $n$ and $0\le j\le \alpha_1n^{2/3}$. 
Thus  (\ref{avar4}) implies that for some $0\le j\le \alpha_1n^{2/3}$, 
$$\xi^{k_j}_{[nx_j]}(nt)
\le ntg(x_j/t) -  (a_2b^2 \delta_1^2/4)  n^{1/3}(\log n)^{2/3}.$$
Recall  that this event is a consequence of 
(\ref{goalineq}), except on an event of probability less than $3\e_0$,
if $n$ is large enough. 
Set $c_2= a_2b^2 \delta_1^2/4$. 
We can  summarize the development
in the inequality 
\beas
&&P\{ X_n(nt)\le nr -bu_n\}\le 3\e_0 \\
&&\qquad\qquad + 
P\{ \xi^{k_j}_{[nx_j]}(nt)
\le  ntg(x_j/t) - c_2 n^{1/3}(\log n)^{2/3} \\
&&\qquad\qquad\qquad\qquad \mbox{ for some $0\le j\le \alpha_1n^{2/3}$}
\}.
\eeas
Apply  Prop.\ \ref{ldboundprop} to the last probability above
with $h=(c_2/2)n^{-2/3}(\log n)^{2/3}$. This probability
 is less
than $\e_0$,  if $b$ (and hence $c_2$) 
is fixed large enough, and if $n\ge n_0$ for some $n_0$ that 
depends on $c_2$. The cutoff $n_0$ is required to make
$h$ small enough for   Prop.\ \ref{ldboundprop}, after 
$c_2$ has been first fixed large enough. The values $n<n_0$ 
can then be accounted for by increasing $b$ suitably.  
This completes the proof of (\ref{upperest1}).

\bibliographystyle{plain}

\end{document}